\documentclass[11pt]{article}
\usepackage{color}
\usepackage{amsfonts}
\usepackage{amsmath}
\usepackage{amsthm}
\usepackage{amscd}
\usepackage{amssymb}
\usepackage[dvips]{graphicx}
\usepackage[all]{xy}
\usepackage{epsfig}
\usepackage{geometry}
\usepackage{setspace}
\usepackage{fancyheadings}
\usepackage{delarray}
\def\circit#1#2#3{\vbox to 0pt{\kern-#2pt{\hbox to 
0pt{\kern#1pt{$_\circ$}\hss}\vss}}#3}

\begin{document}
\newtheorem{lemma}{Lemma}
\newtheorem{defi}{Definition}
\newtheorem{prop}{Proposition}
\newtheorem{cor}{Corollary}
\newtheorem{theorema}{Theorem}
\newtheorem{claim}{Claim}
\newtheorem{fact}{Fact}
\newtheorem{remark}{Remark}
\newtheorem{recall}{Recall}
\newtheorem{example}{Example}

\title{$\mathbb{Z}_n$-manifolds in $4$-dimensional graph-manifolds}
\date{}
\author{A.Mozgova}
\maketitle

\begin{center}
{\tiny Laboratoire d'analyse non lin\'eaire et g\'eom\'etrie,
Universit\'e d'Avignon, 33, rue Louis Pasteur,
84000 Avignon, France.

\texttt{mozgova@univ-avignon.fr} }

\end{center}

\begin{abstract}

A standard fact about two incompressible surfaces in an irreducible
$3$-manifold is that one can move one of them by isotopy so that
their intersection becomes $\pi_1$-injective. 
By extending it on the maps of some $3$-dimensional 
$\mathbb{Z}_n$-manifolds into $4$-manifolds, we prove that any 
homotopy equivalence of 
$4$-dimensional graph-manifolds with reduced graph-structures 
is homotopic to a diffeomorphism preserving the structures. 

\vspace{0.1cm}


\end{abstract}



\section{$\boldmath\mathbb{Z}_n$-submanifolds of $4$-manifolds}

A {\it $\mathbb{Z}_n$-manifold} of dimension $k$ is 
an object $Z$ obtained from a compact oriented $k$-manifold with boundary
by identifying (with respect of the orientation) $n$ isomorphic parts 
of boundary in such a way that if a
point of the boundary component participate in the identification, then
all the points of this component do. 
Every point of $Z$ has a neighborhood isomorphic either to a
$k$-cell, or to a $k$-semi-cell, or to a $n$-page book of $k$-cells. 
The identified parts form {\it the singular set} of $Z$ denoted $Z_s$, 
and the closure of their complement is its {\it regular set} denoted
$Z'$. Non-identified boundary components form the boundary $ \partial Z$.

For example, by identifying boundaries of three oriented surfaces 
each one having one boundary component, one obtain a $2$-dimensional
$\mathbb{Z}_3$-manifold without boundary. 
Mapping cylinder of $n$-fold covering $S^1 \to S^1$ gives an example 
of $2$-dimensional $\mathbb{Z}_n$-manifold with boundary; it appears 
in $3$-dimensional Seifert manifolds 
as the preimage of an arc in the base orbifold going from the projection 
of a singular fiber to the boundary. 

A standard fact about two incompressible surfaces in an irreducible
$3$-manifold is that one can move one of them by isotopy so that
their intersection becomes $\pi_1$-injective. 
Under natural homotopic assumption, this remains true for the map of
$3$-manifold and $3$-dimensional submanifold of $4$-manifold 
(Proposition~2.B.2 of \cite{Stallings-book} and Proposition~1 of \cite{Mozg}).
The following lemma extends it on the maps of some $3$-dimensional 
$\mathbb{Z}_n$-manifolds into $4$-manifolds.

To fix the notations, let $c : S^1 \to S^1$ be $n$-fold covering.
Its mapping cylinder is a $2$-dimensional $\mathbb{Z}_n$-manifolds 
$Map^n(S^1,S^1)= 
(S^1 \times [0,1]) / \mbox{\footnotesize $(S^1 \times \{ 1 \} ) = 
c(S^1 \times \{ 1 \}$)}.$ \\
Let $Z=S^1 \times Map^n(S^1, S^1)$, it is a $3$-dimensional
$\mathbb{Z}_n$-manifold with boundary $\partial Z=S^1 \times 
\partial Map^n (S^1,S^1)$, singular set $Z_s=S^1 \times
Map^n(S^1,S^1)_s$ and regular set $Z'= S^1 \times
(S^1 \times [0;1])$.

\begin{lemma}\label{main-lemma}
Let $W$ be a compact smooth oriented 4-manifold with
$\pi_2(W)=\pi_3(W)=0$ and $M$ be a compact oriented $\pi_1$-injective 
$3$-submanifold with $\pi_2(M)=0$. 
Let  $f: (Z, \partial Z) \to (W, W \setminus M)$ be a $\pi_1$-injective map.

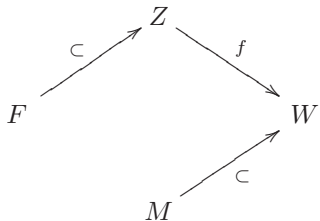
\begin{minipage}[c]{5cm}
{\small
$$\xymatrix{ & Z \ar[rd]^{f} & \\ 
\ \ \ \  F \ \ \ \  \ar[ru]^{\subset} & & \ \ \ \   W \ \ \ \     \\ 
&  M \ar[ru]_{\subset} &   } $$} 
\end{minipage}
\hfill
\begin{minipage}[c]{9cm}
Then one can move $f$ by homotopy (that is constant on $\partial Z$) 
so that $f(Z_s) \cap M = \emptyset$ and each connected component 
of $F=f^{-1}(M)$ is a $\pi_1$-injective torus in $M$.
\end{minipage}
\end{lemma}
\begin{proof}
Move $f$ by a small homotopy to make
it transverse to $M$.
Then $F=f^{-1}(M)$ is a 2-dimensional $\mathbb{Z}_n$-manifold
whose embedding $(F,F_s) \to (Z \setminus \partial Z, Z_s)$ 
induces an embedding of regular sets $F' \to Z'$.
As $M$ is closed and $Z'$ is compact, $F'$ has a finite number 
of connected components (\cite{Abraham-Robbin},
corollary~17.2(IV)).

In $Z = S^1 \times Map^n(S^1,S^1)$ consider the subspace
$G= \{ 0 \} \times Map^n(S^1,S^1)$ where $0 \in S^1$. 
We have $(G, G_s) \subset (Z, Z_s)$ and $G' \subset Z'$.
So $F'$ (which is union of surfaces with boundaries) 
and $G'$ (which is an annulus) are embedded in $Z' = S^1 \times S^1 \times I$. 
Denote the boundary component of $G'$ lying $\partial Z$ by $\partial_0 G'$ 
and the other one by $\partial_1 G'$.

\begin{figure}[ht]
\centering
\input{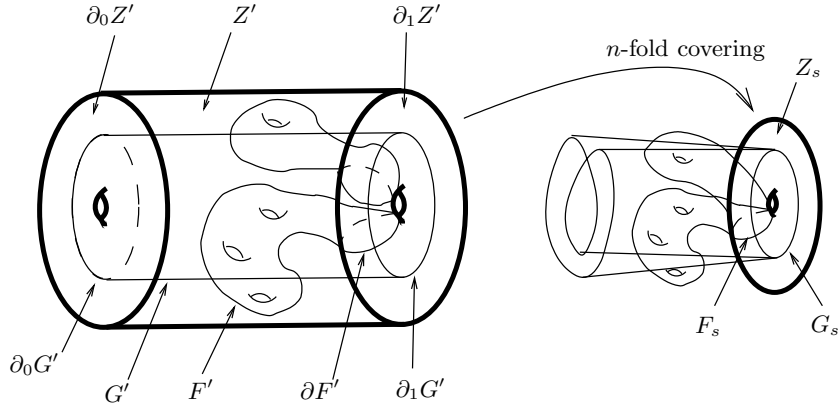}
\protect\caption{{\small Regular parts, the covering $\partial_1 Z' \to
Z_s$, and $Z_s \cup G \cup F$}}
\label{nakr}
\end{figure}

\noindent {\bf Step~1: elimination of trivial circles of $\partial F'$.}
Denote the boundary component of $Z'$ giving $\partial Z$ by $\partial_0 Z'$ 
and the other one (participating in the identification) by $\partial_1 Z'$. 
Suppose that there is a circle of $\partial F'$ that is trivial in
$\partial_1 Z'$. Then its projection on $Z_s$ is trivial, too. As its
projection is embedded in $Z_s \simeq T^2$ (because $\partial F'$ 
is the preimage of $F_s$ which is embedded), it bounds an embedded disk
there. Now we can proceed as in Proposition~1 of \cite{Mozg}: take a map
of $2$-disk in $M$ bounding the same loop. As $\pi_2(W)=0$, the
union of these two disks bounds a map of $3$-disk $\alpha : D^3 \to W$, 
which can be separated from $M$. Denote by $N$ a small book neighborhood
of $D^2 \subset Z_s$ in $Z$. The map $f$ will not be changed on
$\overline{Z \setminus N}$; on $D^2$ the homotopy of $f$ will be the
pushing across $\alpha(D^3)$, and $N \setminus D$ will serve to rely 
the new map on $D^2$ and the old map on $\overline{Z \setminus N}$, 
doing it separately on each leaf of the book. 

\begin{figure}[ht]
\centering
\input{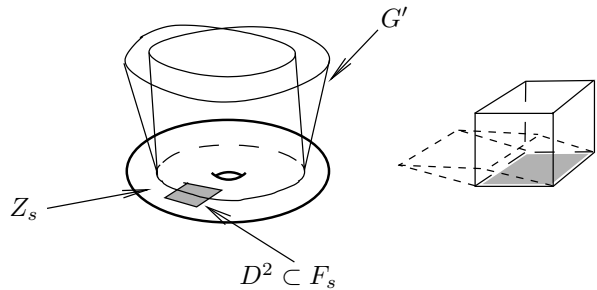}
\protect\caption{{\small Changing of $f$ on $D^2 \subset Z_s$.}}
\label{shape-g}
\end{figure}

\noindent This homotopy of $f$ will eliminate the trivial circle 
from $F'$:

\begin{figure}[ht]
\centering
\input{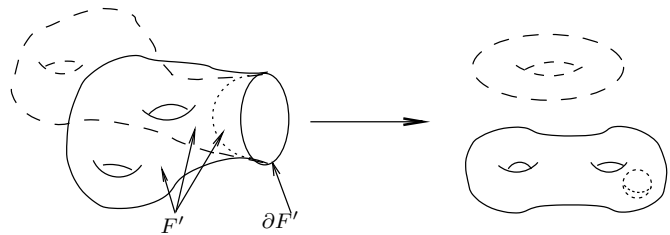}
\protect\caption{{\small Surgery on a trivial loop of $\partial F'$.}}
\label{surg-triv}
\end{figure}

\noindent The closed manifold components of $F$ can be treated 
as in the manifold case.

\noindent {\bf Step~2: reduction of $F'$ to a union of closed 
surfaces and annuli.} 
Now all the circles $\partial F' \subset \partial_1 Z'\simeq T^2$
are non-trivial and embedded, hence, parallel. 
They are either parallel to $\partial_1 G'$, or not, which corresponds
respectively to 
$\partial F' \cap \partial_1 G' = \emptyset$ or 
$\partial F' \cap \partial_1 G' \ne \emptyset$. Take one connected
component of $F'$, denote it still by $F'$.

We will separate two cases: $F' \cap G' = \emptyset$ or 
$F' \cap G' \ne \emptyset$.

\noindent \underline{Case~1: $F' \cap G' = \emptyset$.} 
It implies $F' \subset \overline{\big( Z' \setminus G' \big)} \simeq 
S^1 \times I \times I$. Note that $\overline{ \big( Z' \setminus G'
\big)} \subset Z'$ 
is $\pi_1$-injective but $F' \subset Z'$ is not, so
neither is $F' \subset \overline{ \big( Z'\setminus G' \big)}$.

In this case $(F',\partial F') \rightarrow (S^1\times I\times
I,\partial (S^1 \times I \times I))$, 
$\partial F'$ going to a curve parallel to the generator of
$\pi_1(S^1\times I\times I)=H_1(S^1\times I\times I)$.
We need to find an embedded disk in $S^1\times I\times I$ for the 
trivializing loop.
As $\partial F'$ bounds $F' \subset S^1 \times I \times I$, 
$\partial F' \sim 0$ in $H^1(S^1 \times I \times I,\mathbb{Z}_2)$. 
As all the curves of $\partial F'$
are $\neq 0$ and parallel in $\partial (S^1 \times I \times I)$, they are 
in even number: say $m$ with one orientation, $m$ with the opposite one.
$(F',\partial F') \rightarrow (S^1 \times I \times I,\partial (S^1 \times
I \times I ))$, 
so each $\partial F'$ is sent to a generator of $\pi_1(S^1 \times I
\times I)$.
Take an embedded disk $D = \{ 0 \} \times I \times I \subset 
(S^1 \times I \times I)$. We have then $(D , \partial D) \to 
(S^1 \times I \times I, \partial (S^1 \times I \times I))$ and 
$D \cap F'$ is the union of circles and annuli. 
Take the innermost circle of $D \cap F'$; if it is $\sim 0$ in $F'$, 
we have a disk for the surgery, otherwise move $D$ by isotopy 
to eliminate it, and so on.
After treating all the circles either we have a disk
for the surgery, or there are only arcs in $D \cap F'$. 
Take two components of $\partial F$ corresponding to the arc whose
ends are neighbouring in $\partial D$, denote them by $a$ and $a^{-1}$. 
Then $F'$ car be presented as 
$F'=A \, \raisebox{-1.6mm}{$\stackrel{\bigcup}{\mbox{\tiny $\gamma$}}$} \, F''$ 
where $A$ is an annulus with a hole, $\partial A=a \cup a^{-1} \cup \gamma$
and $F''$ is the remaining part.

\begin{figure}[ht]
\centering
\input{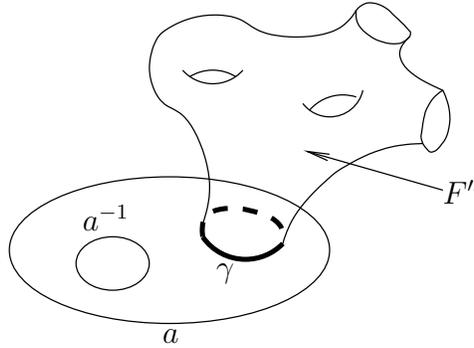}
\protect\caption{{\small $F'$ as union of an annulus with a hole and $F''$.}}
\label{shape-f-pr}
\end{figure}

\noindent Denote by $A'$ the annulus in $\partial (S^1 \times I \times I)$ 
lying between $a$ and $a^{-1}$.
Note that $\gamma \sim aa^{-1} \sim 0$ in $\pi_1(S^1 \times I \times I)$ 
(but $\gamma \nsim 0$ in $F'$ as $F'\neq S^1\times I$).
The arcs $A \cap D$ and $A' \cap D$ bound a disk in $D$, whose
interior does not intersect $F'$.

\begin{figure}[ht]
\centering
\input{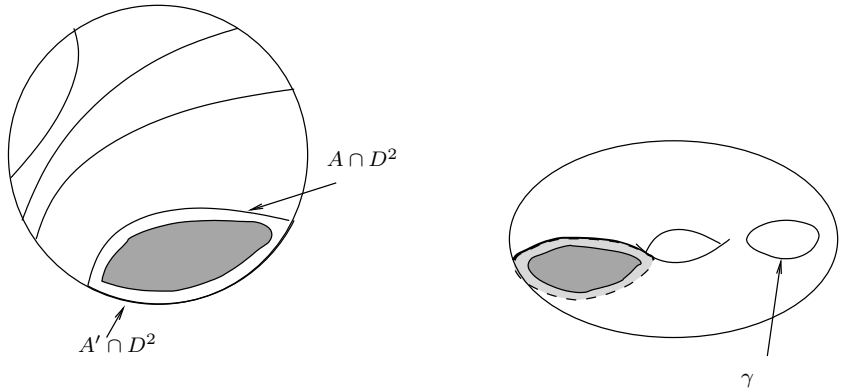}
\protect\caption{{\small $D \cap F'$ in $D$ and $D$ in $S^ 1 \times I
\times I$.}}
\label{a-a-shtr}
\end{figure}

\noindent Push $A \cup A'$ along its normal bundle 
toward the inside of this disk,
living $\gamma$ unchanged, and do the surgery of the pushed part (which
is a torus with a hole) on the
interior disk. We obtain a disk $\tilde{D}$ embedded in $S^1 \times
I \times I$ such that $\partial \tilde{D} = \gamma$ and 
$\tilde{D} \cap F'=\gamma$.
Hence we can move $f$ by homotopy in order to do the surgery of $F'$
on this disk; after this $F'$ becomes the disjoint union of 
an annulus and 
$F'' \, \raisebox{-1.6mm}{$\stackrel{\bigcup}{\mbox{\tiny $\gamma$}}$} \, D^2$.
Doing this for every pair of components of $\partial F'$ corresponding to
an arc with neighboring (in $\partial D$) ends, we will reduce $F'$ to a
union of closed surfaces (which will be treated as in the closed manifold
case) and annuli.

\noindent \underline{Case~2: $F' \cap G' \neq \emptyset$} is a 
$1$-dim manifold (embedded
in both $F'$ and $G'$); recall that  $G'$ was chosen such that
$\big( F_s \cap G_s = \emptyset \Leftrightarrow \partial F' 
\cap \partial_1 G' = \emptyset \big)$.

\noindent\underline{Case~2.1: $F' \cap G' \neq \emptyset$ but 
$F' \cap \partial_1 G' = \emptyset$.}
Then $\partial F'$ is a union of circles parallel to $\partial_1 G'$ ;
and $F' \cap G'$ is a closed $1$-manifold.
Circles of $F' \cap G'$ that are trivial in $\pi_1(F')$, 
can be eliminated by moving $G'$ (without moving $\partial_1 G'$).
So, if all the circles of $F' \cap G'$ are trivial in $\pi_1(F')$, 
we are back to the Case~1.
If there is a circle of $F' \cap G'$ that is non-trivial in $\pi_1(F')$, 
it is either parallel to a component of $\partial F'$ or not.

\noindent $\bullet$ If there exists a circle not parallel to a component
of $\partial F'$, then it is embedded in $G'$ and not parallel to its
boundary, si it bounds a disk in $G'$. We can do the surgery on this disk
(choosing the innermost one).

\noindent $\bullet$ If all circles of $F' \cap G'$ are parallel to a
component of $\partial F'$, then cutting $Z'$ along $G'$, we'll have some
number of $\pi_1$-injective annuli embedded in $S^1 \times I \times I$, 
boundary going to boundary, hence isotopic to annuli in 
$\partial(S^1\times I \times I)$ by an isotopy that is trivial on the
boundary. Use these isotopies in $Z'$ to move $G'$
in order to separate it from $F'$, and we are again in the settings of
the Case~1.

\noindent\underline{Case~2.2:  $F' \cap G' \neq \emptyset$ and 
$\partial F' \cap \partial_1 G' \neq \emptyset$.}
If the components of $\partial F'$ are homotopic to $\partial_1 G'$, move
$G'$ by isotopy to separate $\partial_1 G'$ from $\partial F'$, and we are
in the previous case.
If the curves of $\partial F'$ are not homotopic to $\partial_1 G'$, choose
an embedded curve $\alpha \subset \partial_1 Z'$, parallel to the curves
$\partial F'$ (different from them) and make
$\alpha \times I \subset Z'$ be the new $G'$. Then, we are 
back to Case~2.1.

\noindent {\bf Step~3: the image of the singular set $\boldmath Z_s$ 
can be separated from $\boldmath M$.}
Now every connected component of $F'$ is a torus or an annulus.
Let us show that all the annuli can be eliminated. To simplify the
notations, $F$ will stand for $F$ without tori-components.

Fix a decomposition of $Map^n(S^1,S^1)$ by a wedge of intervals as
follows: 

\begin{figure}[ht]
\centering
\input{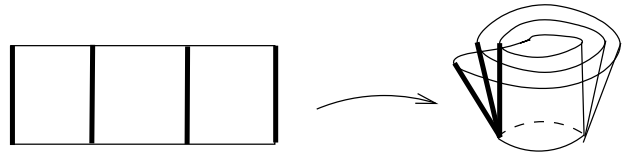}
\protect\caption{{\small Decomposition of $Map^n(S^1,S^1)$ on sheets.}}
\label{decomp}
\end{figure}

\noindent
After multiplying by $S^1$, it gives a wedge of annuli $\bigvee S_i$ 
(identified along one boundary component) embedded in $Z$ 
and decomposing $Z$ into sheets $I \times I \times S^1$, 
in which the corresponding parts of $F$ are manifolds.
One generator of $\pi_1(Z_s)$ is given by the singular circle of 
$Map^n (S^1, S^1)$. The loop $ \big( \bigvee S_i \big) \cap Z_s$ 
corresponds to the second generator, and can be choosen not to be
parallel to $F_s$.
Then every connected component of $F$ intersects a decomposing 
annulus $S_i$ by
arcs. Take an arc coming from a component $F_1 \subset F'$, and suppose 
it is the innermost one in $S_i$ (i.e. its union with an arc from 
$\partial S_i$ bounds a disk in $S_i$ that does not contain other 
arcs from $F \cap S_i$); denote it by $\alpha \subset F_1 \cap S_1$. 
Then the arcs $F_1 \cap S_j$ are the innermost ones 
in $S_j$ $\forall j \ne i$. 
The product of $I$ with the union of the
corresponding disks in $S_1$ and $S_2$ (see Fig.~\ref{sheet-disc} below)

\begin{figure}[ht]
\centering
\input{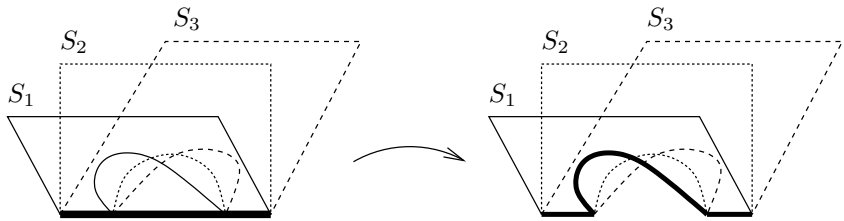}
\protect\caption{{\small Disks in the sheets.}}
\label{sheet-disc}
\end{figure}

\noindent allows to change $f$ 
by homotopy so that the new image of the singular set $Z_s$ coincide 
with the image of $I \times (F_1 \cap S_1)$. Hence by pushing off along
the normal bundle of $M$, $\alpha \times I$ can be eliminated from 
$F \cap \big( \bigvee S_i \big) $ (or, equivalently, the image 
of the sheet of $Z
\setminus Z_s$ that lies beetween $S_1$ and $S_2$ can be pushed 
away from $M$). Taking the product of the union of disks in 
$S_k$ and $S_{k+1}$ by $I$ and doing the same homotopy for the sheet
lying beetween $S_k$ and $S_{k+1}$, one will eliminate all the annuli 
of $F$.

\end{proof}

\section{Application to $4$-dim graph-manifolds}

\subsection{Seifert manifolds}

Following \cite{Ue1, Ue2}, we say that an orientable $4$-manifold $S$ 
with boundary is a {\it Seifert manifold} if
it has the structure of a fibered orbifold $\pi:S \rightarrow B$ 
over a $2$-orbifold with generic fiber $T^2$, and $S$ and $\partial S$ 
are non-singular as orbifolds. Note that $B$ as orbifold has no boundary,
but the underlying surface of $B$ does. 

\noindent {\bf Local picture.}

\noindent Any point $b\in B$ has a neighborhood of type $D=D^2/G$ 
where $G$ is a finite
subgroup of ${\cal O}(2)$ corresponding to the stabilizer of $p$.
Then $\pi^{-1}(D)=(T^2\times D^2)/G$ where the action of $G$ is free
and is a lifting of an action of $G$ on $D^2$
(that is $\pi_{|_{\pi^{-1}(D)}}:(T^2\times D^2)/G \rightarrow D^2/G$
is induced by the canonical projection $T^2\times D^2 \rightarrow D^2$).

\noindent\underline{Case 0}: $G=\{1\}$. \\
In that case $p$ is non-singular and $\pi^{-1}(p)=T^2$ is a regular fiber.
The preimage of an arc in $B$ that joins 
$p$ with $\partial B$ is $T^2\times I$.

\noindent\underline{Case 1}: 
$G=\mathbb{Z}_m= \langle \; g \; | \; g^m=1 \; \rangle$ and the action is
given by $g (x,y,z)=\big(x-a/m,y-b/m,z{\rm e}^{\frac{2\pi i}{m}}\big)$,
where $(x,y)\in T^2$, $z\in D^2=\{z\in\mathbb{C} \; | \; |z| \le 1\}$, 
and $m,a,b$ being mutually relatively prime.
Then $p$ is a cone point of angle $2\pi/m$ and $\pi^{-1}(p)=T^2$.
The covering (regular fiber) $\rightarrow$ (singular fiber)
corresponds to the subgroup 
$a\mathbb{Z}\oplus b\mathbb{Z} \rightarrow \mathbb{Z}\oplus\mathbb{Z}$. 
As $(a,b)=1$, we can choose a basis in the $\pi_1$ of the regular 
fiber such that the covering corresponds to the
subgroup $\mathbb{Z}\oplus n \mathbb{Z}\rightarrow\mathbb{Z}\oplus\mathbb{Z}$.
Then the covering $T^2\times D \rightarrow (T^2\times D^2)/G$
splits as
$$\begin{array}{rcl}
S^1 & \times & (S^1 \times D^2) \\
\mbox{{\tiny id}} \downarrow \hspace*{1.9mm} & & \hspace*{7.1mm} \downarrow \mbox{{\tiny n-fold covering}} \\
S^1 & \tilde{\times} & (S^1 \times D^2) \\
\end{array}$$ 
(the bundle $S^1 \tilde{\times} (S^1\times D^2)$ beeing trivial because
$(T^2\times D^2)/G$ is orientable).
That means that 
$(T^2\times D^2)/G=S^1 \times (\mbox{{\small 3-dimensional model}})$.
So in this case, the preimage of an arc in $B$ that joins 
$p$ with $\partial B$ is
the $\mathbb{Z}_n$-manifold $S^1 \times  Map^n(S^1,S^1)$.

\noindent\underline{Case 2}:
$G=\mathbb{Z}_2= \langle \; \tau \; | \; \tau^2=1 \; \rangle$ and the
action is given by $\tau \cdot (x,y,z)=(x+1/2,-y,\bar{z})$.
Then $\pi^{-1}(p)=K^2$, $p$ lies on a reflector circle and $\pi^{-1}(D)$
is a twisted $D^2$-bundle over $K^2$.
The preimage of an arc in $B$ that joins the point
$p$ with $\partial B$ is in this case a manifold $K^2 \tilde{\times} I$.

\noindent\underline{Case 3}:
$G=D_{2m}= \langle \; \tau , g \; | \; \tau^2=g^m=1 \, , \;
\tau g \tau^{-1}=g^{-1} \; \rangle$ and the action is given by
$\tau \cdot (x,y,z)=(x+1/2,-y,\bar{z})$ and
$g \cdot (x,y,z)=\big(x,y-b/m,z{\rm e}^{\frac{2\pi i}{m}}\big)$.
Then $p$ is a corner reflector of angle $\pi/m$ and
$\pi^{-1}(p)=K^2$ which is $m$-fold covered by $K^2$ over a regular point of 
the reflector circle.

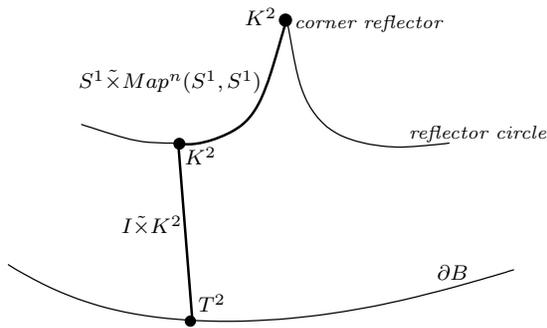
\begin{minipage}{7.5cm}
\vspace{0.5cm}
\input{diedr.pstex_t}
\vspace{0.5cm}
\end{minipage}
\begin{minipage}{8cm}
The preimage of an arc in $B$ that joins $p$ with $\partial B$ is
the mapping cylinder of a covering $T^2 \to K^2$ which is
a 2-fold covering over the loop reversing the orientation of $K^2$
and a $m$-fold covering over the other loop.
This is not a $\mathbb{Z}_k$-manifold, but one can join $p$ with
$\partial B$ by $2$ consecutive arcs $\alpha$ and $\beta$, the first one
lying on the reflector circle and not containing another corner
reflector, the second one joining the other end of $\alpha$ with
$\partial B$. Then $\pi^{-1}(\alpha)=S^1 \tilde{\times} Map^n(S^1, S^1)$
and $\pi^{-1}(\beta)=I \tilde{\times} K^2$.
\end{minipage}
\noindent As in dim 3, on the fundamental group level the singular fibers 
create roots of the regular fiber.

\noindent
{\bf Global picture of $S$ with boundary and hyperbolic base.}

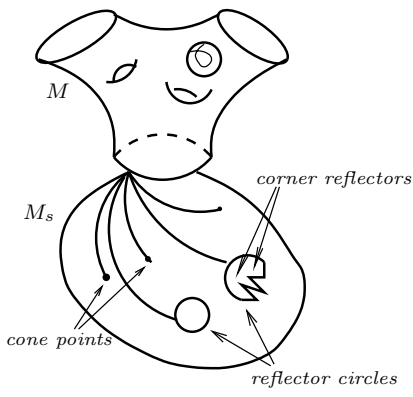
\begin{minipage}{6cm}
\vspace{0.5cm}
\input{glob-pict.pstex_t}
\vspace{0.5cm}
\end{minipage}
\begin{minipage}{9.5cm}
$S=M \cup M_S$, where $M$ is a $T^2$-bundle over a surface and $M_S$ is
its Seifert-part. Note that $0 \rightarrow \pi_1(M) \rightarrow \pi_1(S)$.

\vspace{0.1cm}
\noindent 
{\bf Theorem} (\cite{Zieschang}, \cite{Vogt}, \cite{Ue1}, \cite{Ue2}.)
\textit{Let $S,S'$ be 4-dimensional closed orientable 
Seifert bundles with hyperbolic bases. Then $\pi_1(S)=\pi_1(S')$ 
if and only if there is a fiber-preserving diffeomorphism between $S$ and
$S'$.}

\vspace{0.1cm}
\noindent By taking the doubles, the previous theorem 
implies that if two Seifert manifolds with boundaries 
$(S, \partial S)$ and $(S', \partial S')$ over hyperbolic 
$2$-orbifolds are homotopy equivalent rel boundary 
then $S$ and $S'$ are diffeomorphic.
\end{minipage}

\subsection{$4$-dimensional graph-manifolds}
A block is a Seifert bundle (with boundary) over a hyperbolic $2$-orbifold.
A graph-manifold structure on a compact closed oriented $4$-manifold
is a decomposition as a union of blocks, glued by diffeomorphisms
of the boundary. Note that the boundary of a block has the structure 
of a $T^2$-bundle over a circle.
A graph-manifold structure is
{\it reduced} if none of the glueing maps are isotopic
to fiber-preserving maps of $T^2$-bundles.
Any graph-structure give rise to a reduced one by
forming blocks glued by bundle maps into larger blocks.
Like in the non-singular case, $4$-dimensional graph-manifolds are aspherical,
their Euler characteristic is $0$ (because the blocks are finitely 
covered by $T^2$-bundles over hyperbolic surfaces, 
hence have $\chi=0$, and the glueings are made 
along 3-manifolds), and can be smoothed.

\begin{theorema} 
Any homotopy equivalence
of closed oriented $4$-manifolds with reduced graph-structures 
is homotopic to a diffeomorphism preserving the structures.
\end{theorema}

\noindent We will say that a $\pi_1$-injective map between 
the blocks $f:S \to S'$ 
is {\it fiber covering} if in the fundamental groups it sends 
the (normal) fiber subgroup of $\pi_1(S)$ into the (normal) 
fiber subgroup of $\pi_1(S')$. 
To extend the proof of the non-singular case (\cite{Mozg}), 
one has to show that any $\pi_1$-injective map $f:W =\cup W_i \to W'=\cup W'_k$
of graph-manifolds with reduced graph-structures is homotopic 
to $\bigcup f_i$, where each $f_i : (W_i, \partial W_i) \to (W'_j,
\partial W'_j)$ is a fiber covering map.
The missing step is the following

\begin{prop}
Let $S$ be a block and $W$ be a 4-dim graph-manifold with reduced
graph-structure. Then any $\pi_1$-injective map $f: S \to W$ is homotopic
to a map into one block of $W$. 
\end{prop}

\begin{proof}
In the base of $\pi: S \to B$, take the wedge of circle on which $B$ retracts. 
Join its core by arcs with every cone point and every reflector circle.
As $S$ retracts on 
$$S'=\pi^{-1}\mbox{\small{(wedge of circles + wedge of arcs + 
reflector circles)}},$$ 
\noindent we have to show that $f \arrowvert_{S'}$ is
homotopic to a map into one bloc of $W$.  

\noindent {\bf Step~1: 
$\boldmath f \arrowvert_{\pi^{-1}\mbox{\small{(wedge of circles)}}}$ is
homotopic to a map in one block of $\boldmath W$.} 
To make the reasonnings as in the
singular case, one only needs to show that up to homotopy rel boundary, 
any $\pi_1$-injective map $g: (T^2
\times I, \partial (T^2 \times I)) \to (S, \partial S)$ is either
fiber-covering or is a map into $\partial S$. Note that $S$ is finitely covered
by $\tilde{S}$ which is a $T^2$-bundle over a hyperbolic surface (with
boundary). As $g$ sends boundary to boundary, $Im \, g_{\ast} \subset 
\pi_1(\tilde{S}) \subset \pi_1(S)$, hence $g$ lifts to a map $\tilde{g} :
(T^2 \times I, \partial ) \to (\tilde{S}, \partial \tilde{S})$, 
to which Lemma~2 (\cite{Mozg}) applies.

\noindent 
{\bf Step~2: $\boldmath f \arrowvert_{\pi^{-1}\mbox{\small{(wedge of arcs + 
reflector circles)}}}$ is homotopic to a map in one block of $\boldmath W$.} 
Denote by $M$ the block of $W$ in which 
$f \big( \pi^{-1}\mbox{(wedge of circles)} \big)$ lies 
and denote $p$ the point of $B$ which is the common core of wedge of circles 
and wedge of arcs.
Then $f \arrowvert_{\pi^{-1}(p)} $ is a covering of the regular fiber of $M$.

\noindent
{\bf Step~2.1: arcs to the cone points.} The preimage of such an arc is 
$Z = S^1 \times Map^n(S^1, S^1)$, and $f \arrowvert_{Z}$ sends $\partial
Z$ into $M$, hence $f (\partial Z)$ does not intersect the decomposing
submanifolds of $W$. Apply Lemma~\ref{main-lemma} to $f \arrowvert_{Z}$
and decomposing submanifolds $\bigsqcup M_{\varphi_i}$, we obtain $f
\arrowvert_{Z}(Z_s) \cap  \big( \bigsqcup M_{\varphi_i} \big) = \emptyset$ 
and every 
component of $(f \arrowvert_{Z})^{-1} \big( \bigsqcup M_{\varphi_i}
\big)$ is a torus $\pi_1$-injectively embedded in $Z'=S^1 \times S^1 \times I$.
Hence 
$(f \arrowvert_{Z})^{-1} \big( \bigsqcup M_{\varphi_i} \big)$ is
$\bigsqcup T^2$ that are parallel to $\partial Z$ and hence are sent 
by $f$ on the covering of regular fiber of $M$. Take the first $T_1^2$
from $\bigsqcup T^2$ counting from $\partial Z$, say it comes from
$M_{\varphi} \subset M \cap M'$. If it exists, denote the next torus from
$\bigsqcup T^2$ by $T_2^2$. The restriction of $f$ on $T^2 \times I
\subset Z$ lying between $T^2_1$ and $T^2_2$ is a $\pi_1$-injective map 
$f_1: (T^2 \times I, \partial) \to (M', \partial M')$. As $f(T^2_i)$ is a
covering of the regular fiber of $M$ and the graph-structure is reduced,
$f(T^2_i)$ is not a covering of the regular fiber of the neighbouring
block $M'$, hence $f(T^2_2) \subset M_{\varphi}$ and $f_1$ is homotopic
rel boundary to a map into $M_{\varphi}$. Continuing like this, 
we change $f \arrowvert_Z$ so that 
$(f \arrowvert_{Z})^{-1} \big( \bigsqcup M_{\varphi_i} \big)$ containes
at most one manifold component. If there are none, we are done.
If there is one such component, we have a $\pi_1$-injective map 
$f:(Z, \partial Z, Z_s) \to (M \cup M', M, M')$.

\begin{figure}[ht]
\centering
\input{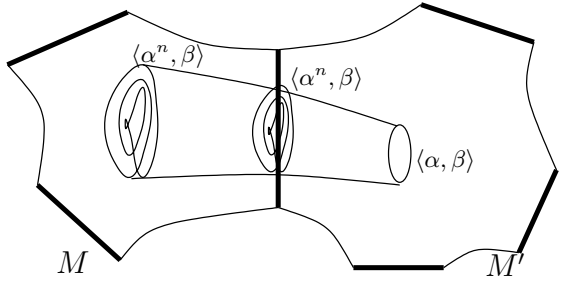}
\protect\caption{{\small $f \big( S^1 \times Map^n(S^1,S^1) \big)$ 
in the blocs $M$ and $M'$}}
\label{blocs}
\end{figure}

\noindent As the torus $(f \arrowvert_Z)^{-1} (\bigsqcup M_{\varphi_i})$ 
is parallel to $\partial Z$, we have a $\pi_1$-injective map 
$f':(Z, \partial Z) \to (M', \partial M')$, such that $f' (\partial Z)$
is not fiber-covering.
  
\noindent \underline{Case~1: $f'$ lifts to a map into $\tilde{M'}$}.
If the torus $\tilde{f'}(\partial Z)$
corresponds to a $(\mathbb{Z} \oplus \mathbb{Z})$-subgroup 
$\langle \alpha^n, \beta \rangle $ in $\pi_1(W)$, then $\tilde{f'}(Z_s)$ 
corresponds to the 
$(\mathbb{Z} \oplus \mathbb{Z})$-subgroup $\langle \alpha, \beta \rangle $. 
Denote $p'$ the projection of $\tilde{M'}$. 
If $p'_{\ast}(\alpha^n)=0$ then 
$p'_{\ast}(\alpha^n)=0$, too because $\pi_1(\tilde{B'})$ is free. Hence
$\alpha$ is homotopic to a loop in the fiber of $\tilde{M'}$.
If $p'_{\ast}(\alpha^n) \ne 0$ then  $p'_{\ast}(\alpha) \ne 0$, and we
can apply 
\begin{quotation}
\noindent {\small 
\begin{fact}
Let $B$ be a surface with boundary, $S$ a component of $\partial B$. 
Then the image $\pi_1(S) \to \pi_1(B)$ is root closed if and only if $B$
is not a M\"{o}bius band and roots are squared.
\end{fact}
\begin{proof}
First note that by \cite{LyndonSch} roots are unique in free
groups. If $\partial B$ has more than one component, then $S$ corresponds 
to a primitive element of $\pi_1(B)$, and the statement comes from
\cite{LyndonSch} ($a^k=s^r$ implies that $a,s$ are powers of a common element,
and as $a$ is primitive, $s$ is a power of $a$).

Now suppose $\partial B$ has one component. Denote $M=Map^k(S^1,S^1)$. 
Suppose there exists $a \in \pi_1(B)$ such that $a \notin 
\pi_1(\partial B)$ and $a^k \in \pi_1(\partial B)$. It gives a
$\pi_1$-injective map $(M, \partial M) \to (B, \partial B)$ such that 
$\pi_1(M) \to a$ and $\pi_1(\partial M) \to a^k$. Attach disks to $M$ and
$B$ to get a $\pi_1$-injective map $M \bigcup_{\partial M} D^2 \to 
B \bigcup_{\partial B} D^2$. As $a^k$ is a finite power of a generator of
$\pi_1(B)$, $\pi_1 \big( M \bigcup_{\partial M} D^2 \big) \to 
\pi_1 \big( B \bigcup_{\partial B} D^2 \big)$ is of finite index.
As $\pi_1 \big( M \bigcup_{\partial M} D^2 \big) = \mathbb{Z}_k$ and 
$\big( B \bigcup_{\partial B} D^2 \big)$ is closed, 
$\big( B \bigcup_{\partial B} D^2 \big) \sim \mathbb{RP}^2$, because  
$\mathbb{RP}^2$ is the only closed surface with finite non-trivial
$\pi_1$. Hence $k=2$ and $B$ is the M\"{o}bius band.
\end{proof} }
\end{quotation}   

\noindent Hence if $p'(\alpha^n)$ is homotopic to a loop in $\partial
\tilde{B'}$, then so is $p'(\alpha)$. Hence $\alpha$ is homotopic to a loop
in $\partial \tilde{M'}$, and $f'$ is homotopic to a map into 
$\partial M'$. 

\noindent \underline{Case~2: $f'$ does not lift to a map into $\tilde{M'}$.}
Then $\alpha^n$ have roots on the singular fibers of $M'$, hence
$\alpha^n$ is homotopic to a loop in the regular fiber of $M'$. 
Denote the Seifert projection of $M'$ by $\pi'$. As the
graph-structure is reduced and $\langle \alpha^n, \beta \rangle$
corresponds to the regular fiber of $M$, 
${\pi}'_{\ast}(\langle \alpha^n, b \rangle) \ne 0$, 
hence ${\pi}'_{\ast}(\beta) \ne 0$. Which means that $\alpha$ can not
commute with $\beta$ because the roots of the regular fiber do not
commute with the elements of the base.

\noindent
{\bf Step~2.2: the image of a reflector circle without corner reflectors
lies in one block}.
The preimage of the reflector circle is twisted $K^2$-bundle over $S^1$.
Applying the previous step to the preimage of the arc between $\partial
B$ and the reflector circle, we see that $f$ restricted to the fiber 
of this bundle is homotopic to a map into the block $M$. Take in $B$ a
loop $\gamma$ near the reflector circle which is homotopic to it.
Then $\pi^{-1}(\gamma)$ is torus bundle $M_{\pm id}$. The image by $f$ of 
its fiber covers the regular fiber of the block $M$, 
hence the whole $f \big( 
M_{\pm id} \big)$ can be shrinked by homotopy into $M$, in particular
$f(\gamma)$ does. Hence the image
of the whole $1$-skeleton of $\pi^{-1}\mbox{(reflector circle)}$ can be
shrinked by homotopy into $M$, and the asphericity implies the statement.  

\noindent
{\bf Step~2.3 : case of reflector circles having corner reflectors.} The
preimage of an arc on a reflector circle between a corner reflector and a
neighboring point is $\tilde{Z}=S^1 \tilde{\times} Map^n(S^1,S^1)$.
Denote as before $\partial \tilde{Z}=S^1 \tilde{\times} \partial \big( 
Map^n(S^1,S^1) \big)$ and $\tilde{Z}_s=S^1 \tilde{\times} \big( 
Map^n(S^1,S^1) \big)_s$. The preimage of an arc in $B$ between the
reflector circle and the boundary is $I \tilde{\times} K^2$, and 
by previous its image  by $f$ can be shrinked into $M$. 

\begin{figure}[ht]
\centering
\input{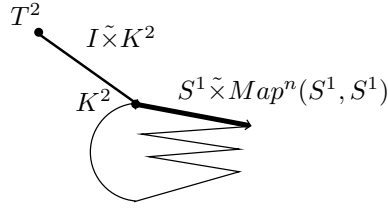}
\protect\caption{{\small Acces to a corner reflector by two arcs.}}
\label{cor-ref}
\end{figure}

\noindent Let us show that
$f(\tilde{Z})$ such that $f(\partial \tilde{Z}) \subset M$ can be
shrinked into $M$. The double covering $p: Z \to \tilde{Z}$ induces a
$\pi_1$-injective map $fp:Z \to W$ which by previous can be shrinked into
$M$ (by a homotopy which is constant on $\partial Z$). Let us show that the
whole mapping cylinder of $p$, $Map_p(Z, \tilde{Z})$ can be shrinked into
$M$. For this, note that the mapping cylinders of the $2$-fold coverings 
$p \arrowvert_{\partial Z}$ and $p \arrowvert_{Z_s}$ are $\big( I
\tilde{\times} K^2 \big)$'s with boundaries lying in $M$, hence both can
be shrinked into $M$. It remains to remark that the union of these
mapping cylinders with $Z$ contains a $1$-skeleton of $Map_p(Z, \tilde{Z})$.

\end{proof}


\begin{thebibliography}{50}


\bibitem{Abraham-Robbin}
R.~Abraham, J.~Robbin, \textit{Transversal mappings and flows},
New York-Amsterdam 1967.-161 p.

\bibitem{Hillman-book}
J.~A.~Hillman, \textit{Four-manifolds, geometries and knots},
Geometry \& Topology Monographs, 5.
Geometry \& Topology Publications, Coventry, 2002.-379 p.


\bibitem{LyndonSch}
R.~C.~Lyndon, M.~P.~Sch\"{u}tzenberger, \textit{The equation
$a\sp{M}=b\sp{N}c\sp{P}$ in a free group}, Michigan Math. J.
\textbf{9} (1962), pp. 289--298.


\bibitem{Mozg}
A.~Mozgova, \textit{Non-singular graph-manifolds of dimension 4}, 
available on \texttt{http://arXiv.org/abs/math.GT/0411335}

\bibitem{Stallings-book}
J.~Stallings, \textit{Group theory and three-dimensional
manifolds}, A James K. Whittemore Lecture in Mathematics given at
Yale University, 1969. Yale Mathematical Monographs, 4. Yale
University Press, New Haven, Conn.-London, 1971.-65 p.

\bibitem{Ue1}
M.~Ue, \textit{Geometric $4$-manifolds in the sense of Thurston 
and Seifert $4$-manifolds. I}, 
J. Math. Soc. Japan 42 (1990), no. 3, pp. 511--540.

\bibitem{Ue2}
M.~Ue, \textit{Geometric $4$ manifolds in the sense of Thurston 
and Seifert $4$ manifolds. II}, 
J. Math. Soc. Japan 43 (1991), no. 1, pp. 149--183.


\bibitem{Vogt}
E.~Vogt, \textit{Vierdimensionale Seifertsche Faserr\"{a}ume},
Dissertation zur Erlangung des Grades eines Doktors der 
Naturwissenschaften der Abteilung f\"{u}r Mathematik an der 
Ruhr-Universit\"{a}t Bochum, Ruhr-Universit\"{a}t, Bochum, 1970. 151 p.


\bibitem{Zieschang}
H.~Zieschang, \textit{On toric fiberings over surfaces}, 
Math. Notes \textbf{5} (1969), pp. 341--345.


\end{thebibliography}
\end{document}